\documentclass[12pt]{article}

\usepackage{amsmath}
\usepackage{graphicx}
\usepackage{amsfonts}
\usepackage{authblk}
\usepackage{url}
\usepackage{float}
\usepackage[numbers]{natbib}

\newcommand{\mycaption}[2]        
 {\begin{center} \parbox{5in}{\caption{\small #2 \label{#1}}} \end{center}}

\newcommand{\pdfgraphic}[4]       
  {\begin{center}
  \scalebox{#2}[#2]{\includegraphics{#1}}
  \end{center}
  \mycaption{#3}{#4}}

\begin{document}

\title{Determining Parameters Leading to Chaotic Dynamics in Systems}

\author[1]{B.C. Dean}
\author[2]{E. Dimitrova}
\author[1]{A. Galande}
\author[2]{E.W. Jenkins}
\author[2]{S. Koshy}

\affil[1]{School of Computing, Clemson University
  (bcdean@clemson.edu)} 
\affil[2]{Department of Mathematical Sciences,
  Clemson University (edimit@clemson.edu,
  lea@clemson.edu, skoshyc@g.clemson.edu)}

\maketitle

\textbf{Abstract:} 
Many biological ecosystems exhibit chaotic behavior, demonstrated either analytically using parameter choices in an associated dynamical systems model or empirically through analysis of experimental data.  In this paper, we provide a computational framework which can be used to both explore the parameter space for the existence of positive Lyapunov exponents and visualize the connections between parameters leading to these positive values.  We demonstrate the effectiveness of the framework on several dynamical systems used to model bacterial populations with a nutrient source.  We provide sample graphics to show the possible ways this framework can be used to gather insight of an underlying system without requiring detailed mathematical analysis. \\ 

\paragraph*{•}

\section{Introduction}
\label{sec:intro}

Although the overall prevalence of chaotic dynamics in biological
systems is not entirely understood, it is well known that chaos occurs
in many natural biological systems involving population dynamics (\citet{costantino97, dennis97, fussmann00, jost73}).  \citet{bec051} demonstrated existence of chaotic and non-chaotic states in an
experimental study of a microbial system with a nutrient source, a
predator species, and two prey species, one of which was less
preferred by the prey .  By measuring the population size
of each species on a daily basis, the authors were able to document
transitions between chaotic and non-chaotic states after changing the
strength of the nutrient source.

Models of population dynamics give us a means of better understanding
chaotic dynamics in biological systems, and even simple dynamical
system models are known to exhibit chaotic behavior for certain values
of system parameters (\citet{kot921, sprott94}).  Unfortunately, even
simple systems can exhibit a range of chaotic and non-chaotic
behavior, related in some non-trivial way to certain combinations of
system parameters and initial conditions.  In this work, we describe a
computational framework that can help researchers develop a better
understanding of the causes of chaotic behavior in dynamical systems.
This platform will allow for efficient systematic exploration of the
parameter space of the system to allow characterization and
visualization of factors leading to chaotic behavior.


We consider models of bacterial populations which incorporate both
predator and prey species, along with a nutrient source.  These
systems have been well-studied in chemostat environments, both
experimentally and analytically
(\citet{banks75,bec051,jost73,kot921}). They have been shown to demonstrate chaotic dynamics
(\citet{bec051, becks08}) and can be considered as small-scale studies
for larger ecological systems(~\citet{banksNotes}).

In general, we consider the set of coupled differential equations
\begin{align*}
\frac{dY_1}{dt} &= f_1 \left( Y_1, \ldots, Y_n, t\right)\\
\vdots\\
\frac{dY_n}{dt} &= f_n \left( Y_1,\ldots, Y_n, t \right)
\end{align*}
where $Y_1 \ldots Y_{n-1}$ represent population levels of various
predator or prey species and $Y_n$ represents the abundance of a
nutrient source.  The functions on the right-hand side are determined
using principles of mass balance, mass action, and enzyme kinetics.
The standard approach uses enzyme-mediated growth rates
(\citet{banksNotes}), where rates are defined as limiting
functions, which cap the ability of the bacterial populations to grow
based on defined saturation-limited values.  These principles are
rooted in Michaelis-Menten kinetics (\citet{banks75,mic131}) and are part
of the standard approaches used to model systems in a chemostat or
other bacterial environment (see, for instance,
\citet{kot921,kov981,str061,hen031,kra041,koo031,molz13} and references
therein).



Given the wide applicability of these models to different physical and
ecological systems, we seek to better understand the choices of model
parameters and initial conditions leading to chaotic behavior.  Let $P
\subseteq \mathbf{R}^m$ denote all possible valid settings for the $m$
parameters in our system (e.g., one parameter might describe the {\rm
  dilution rate} -- the rate at which nutrient is introduced into the
system).  Let $I \subseteq \mathbf{R}^{n}$ denote the set of all valid
initial conditions for $(Y_1, Y_2, \ldots, Y_n)$.  Defining $S = P
\times I$, we seek to characterize the regions of $S$ leading to
chaotic behavior.  A particular $s \in S$ is deemed chaotic if it has
negative divergence at its initial conditions (to ensure boundedness), and yields a positive
Lyapunov exponent.  We use standard numerical integrators for ordinary
differential equations along with popular methods for numerical
Lyapunov exponent calculation (\citet{ben801, wolf85}).

Let $C \subseteq S$ denote the subset of all parameters and initial
conditions leading to chaotic behavior.  The main result of this paper
is a software platform to help researchers characterize the structure
of $C$ for any dynamical system -- a challenging task due to its
high-dimensional nature.  The two main components in our system are
(1) a method for efficiently generating a representative set of values
from $C$ using successive applications of Metropolis-Hastings sampling
to a ``smoothed out'' Lyapunov exponent landscape, and (2) interactive
exploratory visualization of the results on a parallel coordinate
system.  This method of visualization allows the user to observe
correlations among pairs of parameters leading to chaotic behavior,
and also to conduct exploratory testing for specific ``what if''
scenarios -- for example investigating whether chaotic behavior can be
obtained even for highly-restricted ranges of certain parameters.

The remainder of this paper is organized as follows.  In Section
\ref{sec:methods}, we describe the algorithms and visualizations used
in our framework.  In Section \ref{sec:numerics}, we describe the
results of our framework applied to three models for microbial
populations.  Finally, we summarize our findings and provide
discussion on future research directions in Section \ref{sec:conc}.

\section{Computational Methods}
\label{sec:methods}
For a given $s \in S$ describing a set of system parameters and
initial conditions, let $D(s)$ denote the divergence of the system at
$s$ (the trace of its Jacobian evaluated at the initial conditions
found in $s$), and let $L(s)$ denote the maximum Lyapunov exponent of
the system configured according to $s$. If $D(s) < 0$ and $L(s) > 0$,
we say the system exhibits chaotic dynamics.

\subsection{Lyapunov Exponent Calculation}
Lyapunov exponents provide a quantitative measure of the convergence
or divergence of nearby trajectories for a dynamical system. The usual
test for chaos is calculation of the largest Lyapunov exponent: a
positive largest Lyapunov exponent indicates chaos (\citet{sprott03}).
In this work, we use an algorithm developed by  \citet*{wolf85} to compute the entire Lyapunov spectrum
using the solution to the original ordinary differential equation
(ODE) system along with the associated linearized equations of motion.
These latter equations are obtained using the Jacobian of the ODE
system.  Thus, for an $n$-dimensional system, one needs to numerically
resolve a set of $n(n+1)$ differential equations.  The algorithm in
the paper by Wolf, et al., is based on principles found in earlier
works \citet{ben801, shi791}.

One challenge in the use
of numerical integration is the selection of an initial time step and
also a length of time over which to integrate.  Both of these could in
principal be treated as ``system parameters'' and varied automatically
along with all other parameters when searching for chaotic samples.
However, for simplicity, and since the user most likely has some knowledge
of the time scales involved in the system, we ask the user to provide
a time step and also an initial guess for the time range $T$ of
integration.  To make sure the time range is sufficient, we compute
the full Lyapunov spectrum for time range $T, 2T, 4T, 8T, \ldots$,
until the sign of the spectrum does not change between consecutive
doublings (that is, the same number of elements on the spectrum are
positive/negative between doublings, indicating that we have reached a
time range offering at least some level of consistency in our
calculation).

\subsection{Sampling from the Chaotic Regime}
For simplicity, assume for the moment that our system satisfies
$D(s) < 0$ for all $s \in S$, so we want to sample a value $s$
uniformly from the set $C = \{s : L(s) > 0\}$.  We first map $L$
through a {\em sigmoid} function, a continuous approximation of a step
function:
\[
f_{\alpha}(s) = \frac{1}{1 + e^{-\alpha L(s)}}.
\]
The larger we set the parameter $\alpha$, the more $f_{\alpha}$ 
resembles an ideal step:
\[
f_{\infty} = \left\{ \begin{array} {ll}
 1 & \mbox{if } L(s) > 0 \\
 0 & \mbox{if } L(s) < 0. \\
\end{array} \right.
\]
We apply the well-known {\em Metropolis-Hastings} (MH) sampling
algorithm to $f_{\alpha}$.  The MH algorithm performs a random walk
over possible values of $s$, in this case a rectangular region of
$\mathbf{R}^{m+n}$ given by user-specified lower and upper bounds on
every parameter and initial condition.  The MH random walk is designed
so that its stationary distribution is proportional to $f_{\alpha}$,
so after sufficiently many steps, we therefore obtain a sample point
$s$ whose probability of selection is approximately proportional to
$f_{\alpha}(s)$.

During the process of the MH random walk, we slowly increase $\alpha$.
Initially, with $\alpha$ being small, the walk has the flexibility to
move around $S$ in a more fluid manner, so as to avoid getting
``stuck'' in a single local region with $L(s) > 0$ without the ability
to move to other such regions.  As $\alpha$ increases, the walk will
experience greater bias towards staying in regions of $S$ with $L(s)
> 0$, and when $\alpha$ becomes very large the walk will become
essentially trapped within a region of chaotic dynamics, assuming it
found one to begin with.  Since the sigmoid ``flattens out'' the
Lyapunov landscape, our random walk ultimately has the same preference
for choosing any point with $L(s) > 0$, irrespective of the magnitude
of $L(s)$, as just the sign of $L$ is ultimately all we care about when
characterizing chaotic behavior (i.e., a larger positive value of $L$
does not necessarily mean a point is ``more chaotic'').

The method above for generating a single sample is highly reminiscent
of the popular {\em simulated annealing} optimization algorithm, and
indeed we can think of this approach as being an application of
simulated annealing with a goal of maximizing $f_{\infty}$.  This is
not particularly surprising, given the underlying random walk
structure of both the MH and simulated annealing techniques.  If all
we wanted to find was a single sample from the chaotic regime, we
might have characterized our approach more as an optimization-based
method based on simulated annealing.  However, since we want to
generate multiple samples $s \in C$ in order to characterize the
structure of $C$, we feel our approach is more accurately
characterized as one of sampling, based on the MH algorithm.

\begin{figure}[t]
\pdfgraphic{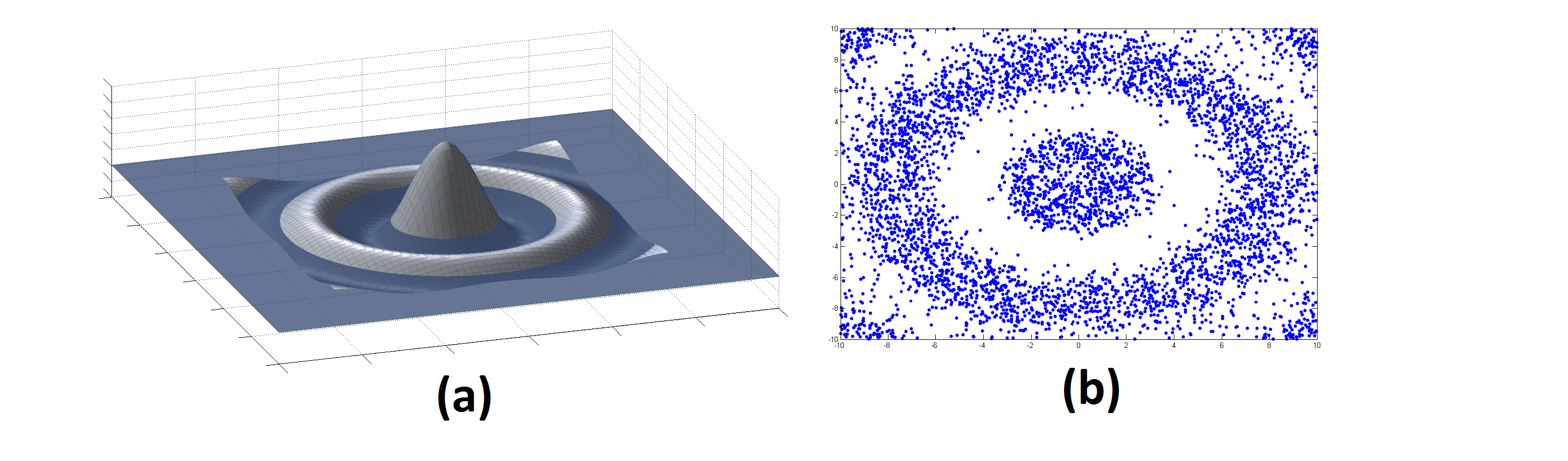}{0.2}{fig_sample}{Sampling 5000 points
from the region on which a test function takes positive values.}
\end{figure}

An example of our approach applied to a test function is shown in
Figure \ref{fig_sample}.  Here, each of 5000 samples is obtained after
taking 1000 steps of an MH random walk, while simultaneously varying
the parameter $\alpha$ from zero up to 20.  Observe that nearly all
samples are drawn from region on which the test function takes
positive values.  

From an efficiency standpoint, each sample drawn by the MH algorithm
involves many steps of a random walk, each of which entails a full
Lyapunov exponent calculation.  Each MH sample generally takes on the
order of several seconds to compute, depending on the complexity of
the underlying dynamical system.  Fortunately, the process of
generating $k$ independent MH samples is ``embarassingly parallel'',
and by parallel implementation on a supercomputer (we use the
``Palmetto'' Cluster at Clemson University), we can obtain all $k$ MH
samples simultaneously in roughly the same time it takes to produce a
single sample.

We now re-introduce divergence into the picture.  Let $C' = \{s : D(s)
< 0\}$ and let $C = \{s \in C' : L(s) > 0\}$.  We now find each sample
point using two ``phases''.  The first samples uniformly a value of $s
\in C'$ using the same smoothed MH random walk, only applied to a
sigmoid computed from $D$ instead of $L$.  If the first phase
succeeds, we proceed to a second phase where we sample from $C$
conditioned on membership to $C'$.  This is done as before, only we
set the probability of a step out of $C'$ to zero during the MH random
walk.  The resulting point is a sample uniformly selected from $C$,
our desired chaotic regime.

\begin{figure}[t]
\pdfgraphic{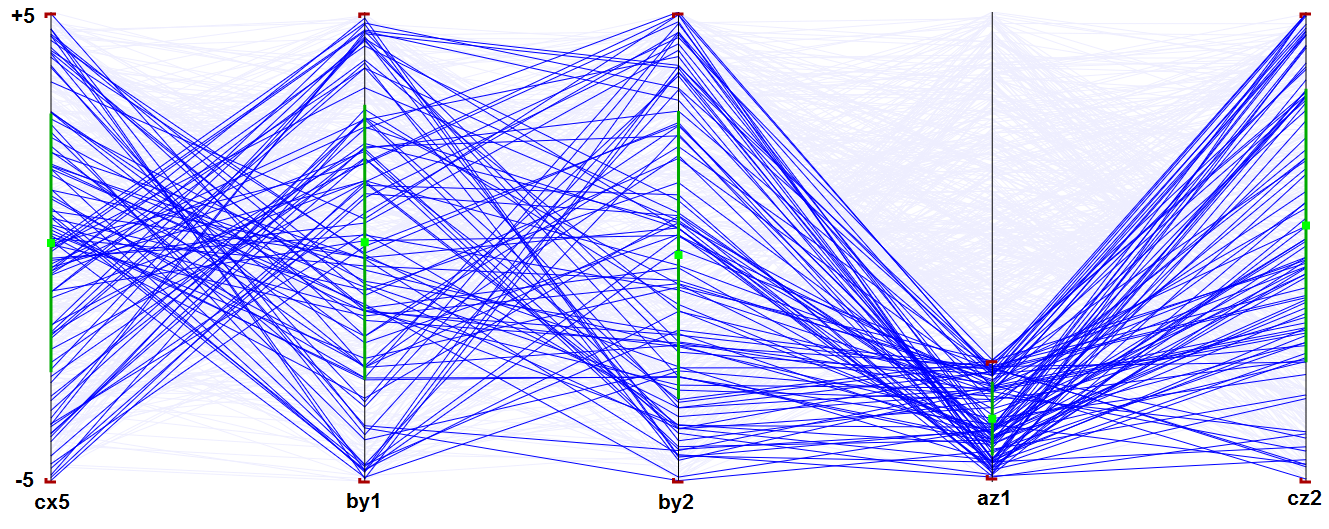}{0.4}{fig_pc}{Chaotic samples visualized
using parallel coordinates.  These results were obtained using one of the systems described in the work by  \citet{sprott94}.  Note the wide range of values for $b_{y1}$ and $b_{y2}$, while the range for $a_{z1}$ is smaller.}
\end{figure}

\subsection{Interactive Visualization}
Now that we have a method for efficiently generating samples from the
chaotic regime $C$, our second main contribution is a proposal for an
interactive visualization platform to help the user understand the
structure of $C$ --- a challenging task, due to the potentially
complex, high-dimensional nature of this set.  Our platform is 
based on {\em parallel coordinates}, a popular means of
visualizing high-dimensional data.  

Figure \ref{fig_pc} shows an example of chaotic sample points
generated from the equations of Case N in Table I of \citet{sprott94}.
The general system in \citet{sprott94} is given as:
\begin{align*}
\frac{dx}{dt} &= a_{x1}+b_{x1}x+b_{x2}y+b_{x3}z+c_{x1}x^{2}+c_{x2}xy+c_{x3}xz+c_{x4}y^{2}+c_{x5}yz+c_{x6}z^{2}\\
\frac{dy}{dt} &= a_{y1}+b_{y1}x+b_{y2}y+b_{y3}z+c_{y1}x^{2}+c_{y2}xy+c_{y3}xz+c_{y4}y^{2}+c_{y5}yz+c_{y6}z^{2}\\
\frac{dz}{dt} &= a_{z1}+b_{z1}x+b_{z2}y+b_{z3}z+c_{z1}x^{2}+c_{z2}xy+c_{z3}xz+c_{z4}y^{2}+c_{z5}yz+c_{z6}z^{2}
\end{align*}
    
The 5 parameters we varied are $c_{x5}$, $b_{y1}$, $b_{y2}$, $a_{z1}$
and $c_{z2}$; the remainder of the parameters were chosen to be zero. Each of the selected varied parameters
is mapped to an individual coordinate axis, arranged in
parallel from left to right, with each sample point drawn as a ``poly
line'' that intersects each coordinate axis at the appropriate
location.  All five parameters were constrained to lie in
the interval $[-5,5]$ during initial sampling, although by
dynamically dragging the upper and lower endpoint markers on each
axis, the user can further restrict the display so it only shows
samples generated within a smaller sub-rectangle.  Here, we have
decreased the upper bound on the parameter $a_{z1}$ to $-2.4$, filtering out
some of the $500$ initial sample points initially present.

\begin{figure}[t]
\pdfgraphic{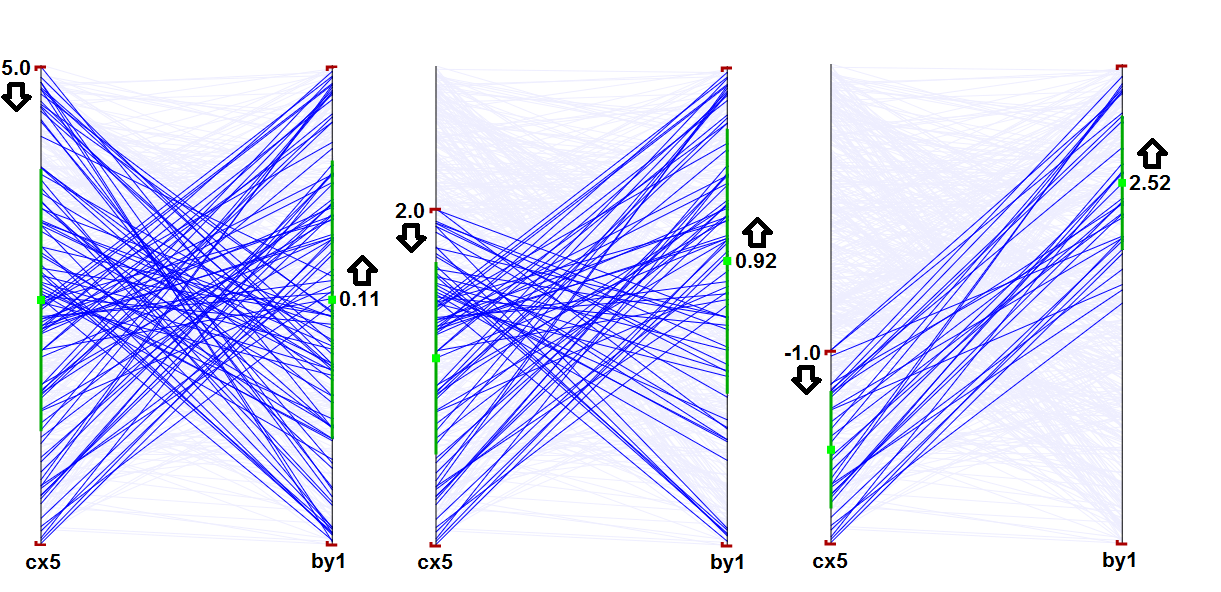}{0.4}{fig_drag}{The result (from left to right) of dragging 
down the upper limit on the $c_{x5}$ axis, observing a corresponding increase
in the mean of the $b_{y1}$ axis due to anticorrelation between the two.}
\end{figure}

Parallel coordinate plots allow us to understand a number of useful
properties by visual inspection and interactive manipulation.  For
example, in Figure \ref{fig_pc}, we see that in order to achieve chaos
with $a_{z1}$ restricted to such a low range, the parameters $c_{x5}$
and $b_{y1}$ need to be anticorrelated: that is, there needs to be a
negative relationship in which $c_{x5}$ increases as $b_{y1}$
decreases. This is visually apparent from the ``X'' pattern between
the $c_{x5}$ and $b_{y1}$ axes, and we can also see it by dynamically
dragging the upper bound on $c_{x5}$ downward, watching the marker for
the mean value on the $b_{y1}$ axis move upward in lock step, as shown
in Figure \ref{fig_drag}.  Interestingly, if one raises the lower
bound on $a_{z1}$, forcing this parameter to take large values, then
the pattern between $c_{x5}$ and $b_{y1}$ becomes one of mostly
straight horizontal lines, indicating correlation rather than
anticorrelation.

By restricting several coordinates at a time, the user can filter an
initially large number of sample points down to only a few.  For
example, if we restrict the range to half of each of five axes, this
will on average show only $1/2^5$ of all our sample points.  In order
to populate the filtered space with sufficiently many samples to
understand its geometric structure, it may be necessary to re-launch a
new round of MH sampling within this restricted space, requiring tight
coupling between the user interface and the ``back end'' parallel MH
sampler.

\section{Numerical Results}
\label{sec:numerics}

As the main goal of our work is to provide researchers with a
convenient numerical method for exploring chaotic behavior of their
systems of interest, we tested the ability of our algorithm to find
parameter settings leading to chaotic states using three systems from
literature.  We use this section to describe those systems and provide
the associated numerical results.  At the end of every example we make
observations on the combinations of parameters and initial values that
lead to chaos as examples of the type of information that researchers
studying such systems will be able to derive using our methods without
rigorous mathematical study.

We begin with the system presented by  \citet{kot921},
where the authors analyze a forced double-Monod model in an initial
effort to understand chaos in biological systems.  We next consider a
system of equations for modeling plant growth in the rhizosphere
(\citet{kra041, str061}), which includes multiple prey species, a
predator species, and a nutrient source.  Finally, we use our
framework to analyze a system of equations motivated by experimental
measurements of a bacterial population with a nutrient source
(\citet{bec051}).


\subsection{Double forced Monod system}

We first consider one of the dynamical systems presented and analyzed
in work by  \citet*{kot921}.  While experimental
results are not included in the paper, the authors justify their
choice to study this particular system by noting the possibility of
obtaining experimental validation of their work.  As our primary
interest in developing the framework presented here is to aid
biologists in their data-based studies of dynamical systems, this
particular problem presents an ideal benchmark case.

The dimensionalized equations are given as 
\begin{align*}
\frac{dS}{dt} &= D\left[ S_i \left(1 + \epsilon
  \sin\left(\frac{2\pi}{T}t\right)\right) - S\right] -
\frac{\mu_1}{Y_1} \frac{SH}{K_1+S} \\
\frac{dH}{dt} &= \mu_1 \frac{SH}{K_1+S} - DH - \frac{\mu_2}{Y_2}\frac{HP}{K_2+H} \\
\frac{dP}{dt} &= \mu_2 \frac{HP}{K_2+H} - DP,
\end{align*}
where $S$ represents the limiting substrate, $H$ represents a prey
species, and $P$ represents a predator species.  We note the predator
species consumes only the prey, so its population is indirectly
associated with the changes in $S$.  The parameters in the model
govern the response of the organisms to changes in the system.  $D$ is
the dilution rate, which defines the ratio of the flow into the
chemostat to the volume of organisms in the chemostat; $\mu_1$ and
$\mu_2$ are the maximum specific growth rates of the prey and
predator, respectively; $Y_1$ is the yield of prey per unit mass of
substrate, and $Y_2$ is the yield of predator per unit mass of prey;
and $K_1$ and $K_2$ are the half-saturation constants.
$S_i$ is the inflowing substrate concentration.

The equations are nondimensionalized by rescaling all variables by
the inflow substrate, the prey by its yield constant $Y_1$, and the
predator by both yield constants (\citet{kot921}).  The dimensionless
system is given by
\begin{align*}
\frac{dx}{d\tau} &= 1 + \epsilon \sin\left(\omega \tau\right) - x - \frac{Axy}{a+x} \\
\frac{dy}{d\tau} &= \frac{Axy}{a+x} - y - \frac{Byz}{b+y} \\
\frac{dz}{d\tau} &= \frac{Byz}{b+y} - z,
\end{align*}
where $x=\displaystyle\frac{S}{S_i}$,
$y=\displaystyle\frac{H}{Y_1S_i}$, and
$z=\displaystyle\frac{P}{Y_1Y_2S_i}$, $\tau=Dt$, and
$\omega=\displaystyle\frac{2\pi}{DT}$.  Also,
$A=\displaystyle\frac{\mu_1}{D}$, $a=\displaystyle\frac{K_1}{S_i}$,
$B=\displaystyle\frac{\mu_2}{D}$, and
$b=\displaystyle\frac{K_2}{Y_1S_i}$.

The dimensionless system exhibits chaotic behavior for
$\omega=\displaystyle\frac{5\pi}{6}$ and $\epsilon=0.6$, with
$S_i=115$ mg/l, $D=0.1$ / h, and initial conditions $x=0.42 \text{, }
y=0.4 \text{, } z=0.42$. The remainder of the parameter values are
provided in Table~\ref{tab:kotParams}. Representative manifold plots
are given in \citet{kot921} and are not duplicated here.

\begin{table}[h!]
\begin{center}
\begin{tabular}{cccc}
\hline
 & $Y_i$ & $\mu_i$ h$^{-1}$ & $K_i$ mg/l \\
 \hline
Prey ($i=1$) & 0.4 & 0.5 & 8 \\
Predator ($i=2$) & 0.6 & 0.2 & 9 \\
\hline
\end{tabular}
\caption{Parameter values associated with the chaotic system described in~\citet{kot921}.}
\label{tab:kotParams}
\end{center}
\end{table}

We wanted to replicate the findings in \citet{kot921} in order to
confirm that our framework can identify parameter values leading to
chaotic and non-chaotic states.  The authors  \citet{kot921} noted
that particular values of $\omega$ and $\epsilon$ were associated with
chaotic states, so we chose to search over ranges of these two
parameters, along with ranges for the initial conditions for the
system.  A subset of our results is shown using the parallel
coordinates visualization in Figure \ref{PCKot}.  This visualization
indicates more positive Lyapunov exponents associated with small
values of $\omega$. Notice that in the rescaled system $\omega$ is the
angular frequency of the forcing term. Observing positive Lyapunov
exponents for small values of $\omega$ may indicate to a biologist an
interesting connection between periodic forcing and chaotic behavior
that they may be able to further study experimentally.  Any clustering
of the remainder of the parameters chosen for our study is less
obvious.  Figure \ref{Manifold_Kot} contains a manifold plot of the
dimensionless system for a randomly chosen set of parameter values
$\epsilon = 0.4809$, $\omega = 1.9109$, and initial conditions $x(0) =
0.3878$, $y(0)=0.7866$, and $z(0)= 0.8875$.  This parameter set
returned a positive Lyapunov exponent.  We note that the strange
attractor is distinct from the one given in the referenced paper.

\begin{figure}[h!]
\pdfgraphic{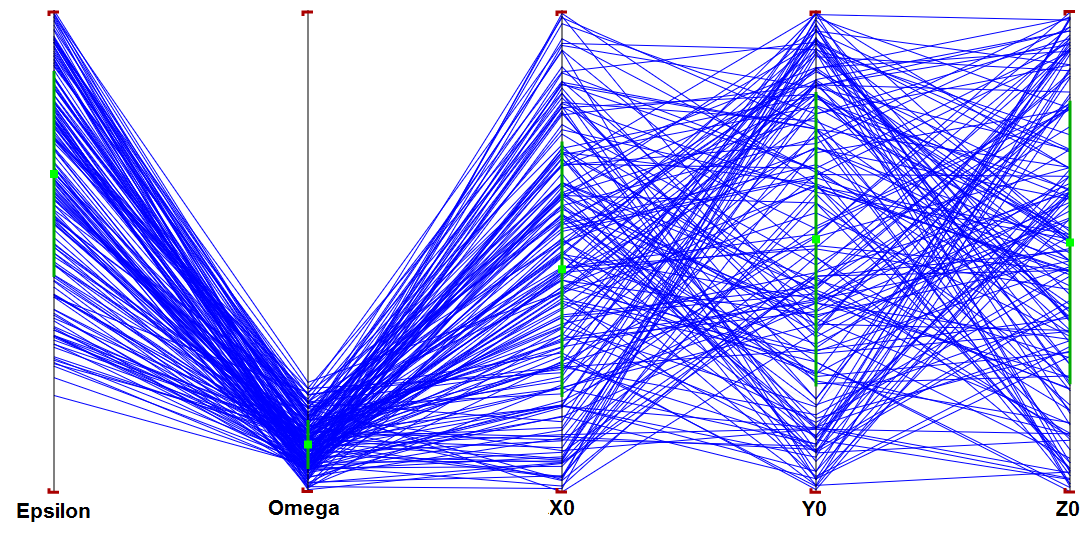}{0.5}{PCKot}{Results from the parameter search
  for the double forced Monod system in \citet{kot921} visualized using
  parallel coordinates.  The vertical axes are associated with values
  for $\epsilon$, $\omega$, $x(0)$, $y(0)$, and $z(0)$, and the
  maximum calculated Lyapunov exponent.}  
\end{figure}

\begin{figure}[h!]
\pdfgraphic{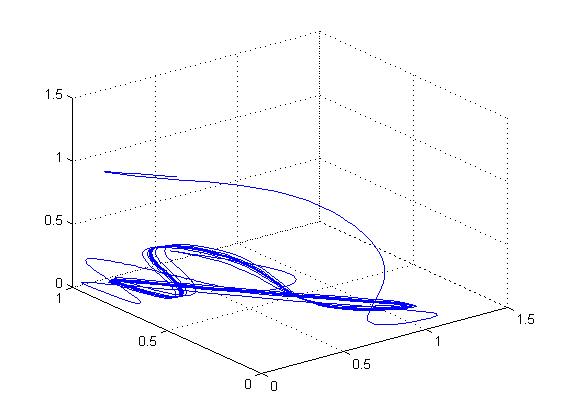}{0.6}{Manifold_Kot}{Manifold plot of
  forced model in \citet{kot921}.  Note the strange attractor is
  distinct from the one given in the paper.}
\end{figure}

\subsection{Models of rhizozone bacteria}

Our second example system is motivated by studies of bacterial growth
in the rhizosphere.  In particular, we consider the models for plant
growth promoting rhizobactera (PGPR) inoculation developed by
 \citet*{kra041} and  \citet*{str061}. In those works, the authors developed the
following model equations to describe the dynamics of their two prey,
one predator system with a periodic nutrient source:
\begin{align*}
\frac{dX}{dt} &= X \left(\mu_X\left[S,P,N\right] + F\left[Z\right] -
\alpha X - d_1\right) \\
\frac{dZ}{dt} &= Z \left(\mu_z\left[S,P,N\right] + G\left[X\right] -
\beta Z - d_2\right) \\
\frac{dS}{dt} &= W\left(t\right) + L - D_S \left(S-S_0\right) -
\frac{X\mu_x\left[S,P,N\right]}{Y_{XS}} -
\frac{Z\mu_z\left[S,P,N\right]}{Y_{ZS}} \\
\frac{dP}{dt} &= D_P \left(P_0 - P\right) -
\frac{X\mu_x\left[S,P,N\right]}{Y_{XP}} -
\frac{Z\mu_z\left[S,P,N\right]}{Y_{ZP}} 
\label{eq:kravchenko}
\end{align*}
In this system, $X$ and $Z$ represent the concentration of PGPR and resident
micro-organisms, respectively, $S$ represents an organic substrate
concentration, and $P$ represents the concentration of oxygen in the
soil (\citet{str061}).  The parameters $\mu_X$ and $\mu_Z$ are growth
rates dependent on the available substrate $S$, oxygen $P$, and
nitrogen $N$ through defined growth rate functions.  These functions
are rate limited and have the general form 
\begin{equation*}
\mu_*\left[S,P,N\right]=\mu_{m*} \frac{S}{S+\theta K_{S*}}
\frac{P}{P+K_{P*}}\frac{N}{N+\theta K_{N*}},
\end{equation*}
where $\mu_{m*}$ is a specified maximal growth rate for the species of
interest, $\theta$ is the moisture content of the soil, and $K_{S*}$
represents an affinity constant for the species to the organic
substrate.  The constants $K_{N*}$ and $K_{P*}$ are similarly defined.
The initial values of the parameters used in our work are provided in detail in~\citet{str061}.

The function $W(t)$ is meant to simulate the effects of
photosynthesis, and is therefore periodic with period length $24$ h.
We incorporate this into our model system using a Fourier series
expansion of the piecewise constant function
\begin{equation*}
f(x)=\left\{\begin{array}{rl}
1, & 0<x\leq 12 \\
0, & 12< x < 24. \end{array}\right.
\end{equation*}

Mathematical analysis of the PGPR models was not provided in the
papers we considered (\citet{kra041, str061}).  However, the authors of
those works indicate that rigorous analysis of the system should prove
the existence of chaotic states.  In addition, as the characteristics
of the PGPR models are similar to other dynamical systems known to
exhibit chaotic behavior, we believe this to be a reasonable system
for further evaluation of our framework.  Towards that end, we use our
algorithm to search over ranges of $K_{SX}$ and ranges of initial
values for $X$, $Z$, $S$, and $P$.  A subset of the results we obtain
is shown in the parallel coordinates plot in Figure~\ref{PCKrav}.  We
note the definite clustering of values for the parameters.  In
particular, we note the value of $Z(0)$ is almost uniformly on the
lower end of the selected search range. Such a finding would
facilitate a researcher interested in studying (or avoiding) chaotic
behavior of their system by clearly pointing to the low values to
which $Z(0)$ should (not) be set.

\begin{figure}[h!]
\pdfgraphic{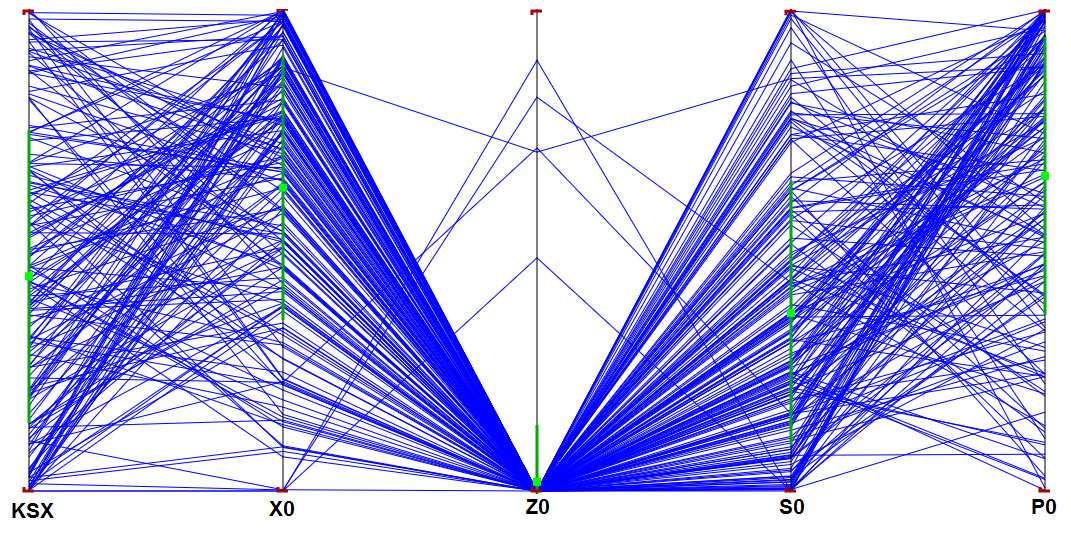}{0.5}{PCKrav}{A parallel coordinates plot for a subset of the results for the the PGPR models in \citet{str061}.  The vertical axes are associated with values of $K_{SX}$, $X(0)$, $Z(0)$, $S(0)$, and $P(0)$, along with the associated maximum Lyapunov exponent.}
\end{figure}

\subsection{Models for experimental data}

For the third example, we model a set of experimental data where
chaotic states were demonstrated for varying levels of the dilution
parameter (\citet{bec051}).  The authors of that work describe
experimental scenarios where the system under study could transition
from chaotic to equilibrium states, and vice versa, by crossing a
threshold value for the dilution rate.

Our approach for modeling the data uses a set of rate-limited,
dynamical equations similar to those in the previous examples. The
model equations and system parameters were estimated by  in~\citet*{molz13}.  Their intent was to derive a
mathematical model whose dynamics closely resembled the experimental
dynamics seen in~\citet{bec051}.  The general description of the system
we consider is given by
\begin{eqnarray}\label{eq:becks}
\begin{aligned}
\frac{dR}{dt} &= R \left[ \mu_{NR} \left( \frac{N}{K_{NR}+N}\right) - \delta_R\right] - \frac{\mu_{PR}}{Y_{PR}} \left(\frac{R}{K_{PR} + R}\right) P - DR \\
\frac{dC}{dt} &= C \left[ \mu_{NC} \left( \frac{N}{K_{NC}+N}\right) - \delta_C \right] - \frac{\mu_{PC}}{Y_{PC}} \left(\frac{C}{K_{PC}+C}\right) P - DC \\
\frac{dP}{dt} &= P \left[ \mu_{PR} \left( \frac{R}{K_{PR}+R}\right) + \mu_{PC} \left(\frac{C}{K_{PC}+C}\right) - \delta_P \right] - DP  \\
\frac{dN}{dt} &= DN_0 - R\left[ \frac{\mu_{NR}}{Y_{NR}} \left(\frac{N}{K_{NR}+N}\right) \right] - C \left[ \frac{\mu_{NC}}{Y_{NC}}\left( \frac{N}{K_{NC}+N} \right) \right] - DN. 
\end{aligned}
\end{eqnarray}
The variables $R$ and $C$ represent the prey species of rods and
cocci, respectively.  We let $P$ represent the prey species, and $N$
represents a nutrient source to the system.  The parameter $D$
represents the dilution rate for input of nutrients to the system.

The parameters in the model determine the feeding habits of the
predator and prey, as well as death and growth rates for each species.
These parameters may be used to specify particular behaviors of the
organisms, e.g., growth rates due to feeding on nutrient sources
rather than prey.  We define $\mu_{N*}$ as maximum growth rates for
the associated species based on consumption of nutrients and
$\mu_{P*}$ as the maximum growth rates for predator based on
consumption of the associated prey species.  The value of $K_{N*}$ is
the half saturation constant for the species on the nutrient, and
$K_{P*}$ is the half saturation constant for the predator on the
associated species. Note that these latter constants may determine a
``preference'' for one prey over the other.  The parameters $Y_{N*}$
represent yield coefficients for the species on the nutrient, and
$Y_{P*}$ represents the yield coefficients for the predator associated
with the prey species.  Death rates for each species are given by
$\delta_*$.  We provide specific values we used for these parameters
in Table~\ref{tab:becksParams}.

\begin{table}[h!]
\begin{center}
\begin{tabular}{|c|c|c|c|}
\hline
 & \multicolumn{3}{|c|}{Species}\\
\hline
Parameter & R & C & P \\
\hline
$\mu_{N*}$ & 12 / day & 6 / day & \\
$\mu_{P*}$ & 2.2 / day & 2.2 / day & \\
$K_{N*}$ & 8e-6 gm/cc & 8e-6 gm/cc & \\
$K_{P*}$ & 1e-6 gm/cc & 1e-6 gm/cc & \\
$Y_{N*}$ & 0.1 gm R / gm N & 0.1 gm C / gm N & \\
$Y_{P*}$ & 0.12 gm P / gm R & 0.12 gm P / gm C & \\
$\delta_*$ & 0.5 / day & 0.25 / day & 0.08 / day \\
\hline
\end{tabular}
\caption{Parameter values for model equations \eqref{eq:becks} used in the numerical simulations.  \label{tab:becksParams}}
\end{center}
\end{table}

The model incorporates fifteen parameters.  The
parameter space can be reduced by rescaling the state variables $R$,
$C$, $P$ and $N$, and time $T$.  Rescaling the equations eliminates
four of the parameters and simultaneously makes the system more
perspicuous for both qualitative analysis and numerical
simulation. We chose to rescale by introducing the following change of variables
\begin{align*}
 R &= K_{PR}\,r, \quad C = K_{PC}\,c, \quad P=\dfrac{K_{PR}Y_{PR}\left(\delta_{R}+D\right)}{\mu_{PR}}\,p, \\ 
 N &= K_{NR}\,n, \quad T = \dfrac{1}{\delta_{R}+D}\,t,
\end{align*} 
where the lowercase $r$, $c$, $p$, $n$ and $t$ are the new rescaled
variables and time, respectively. The rescaled model equations are given by

\begin{eqnarray}\label{eq:snutrient}
\begin{aligned}
\frac{dr}{dt} &= \widehat{\mu}_{NR} \left( \frac{nr}{n+1}\right) -
\left(\frac{rp}{r+1}\right) - r  \\
\frac{dc}{dt} &= \widehat{\mu}_{NC} \left(
\frac{nc}{n+\widehat{\kappa}}\right) - \widehat{\eta}_{1} \left(
\frac{cp}{c+1}\right) - \widehat{\delta}_{C}\, c \\
\frac{dp}{dt} &= \widehat{\mu}_{PR} \left( \frac{rp}{r+1}\right) +
\widehat{\mu}_{PC} \left(\frac{cp}{c+1}\right) -
\widehat{\delta}_P\,p  \\
\frac{dn}{dt} &= \widehat{\delta}\left(n_0 - n\right) -
\widehat{\eta}_{2}\left( \frac{nr}{n+1}\right) - \widehat{\eta}_{3}
\left(\frac{nc}{n+\widehat{\kappa}} \right), 
\end{aligned}
\end{eqnarray}
where the new (eleven) parameters are expressed in terms of the original parameters as follows:
\begin{align*}
\widehat{\mu}_{NR} &= \dfrac{\mu_{NR}}{\left(\delta_{R}+D\right)},
\quad \widehat{\mu}_{NC} =
\dfrac{\mu_{NC}}{\left(\delta_{R}+D\right)}, \\
\widehat{\mu}_{PR} &= \dfrac{\mu_{PR}}{\left(\delta_{R}+D\right)},
\quad \widehat{\mu}_{PC} =
\dfrac{\mu_{PC}}{\left(\delta_{R}+D\right)}, \\
\widehat{\kappa} =\dfrac{K_{NC}}{K_{NR}},\, \widehat{\delta}_{C} &=
\dfrac{\delta_{C}+D}{\left(\delta_{R}+D\right)},\,
\widehat{\delta}_{P} =
\dfrac{\delta_{P}+D}{\left(\delta_{R}+D\right)},\, \widehat{\delta} =
\dfrac{D}{\left(\delta_{R}+D\right)}, \\
\widehat{\eta}_{1} = \dfrac{\mu_{PC}Y_{PR}K_{PR}}{\mu_{PR}Y_{PC}K_{PC}}, \, 
\widehat{\eta}_{2} &= \dfrac{\mu_{NR}K_{PR}}{Y_{NR}K_{NR}\left(\delta_{R}+D\right)}, \, 
\widehat{\eta}_{3} = \dfrac{\mu_{NC}K_{PC}}{Y_{NC}K_{NR}\left(\delta_{R}+D\right)}.
\end{align*}
We also introduce the initial value for the rescaled nutrient source, $n_{0}= K_{NR} N_{0}$. 

We studied the dynamics of our proposed model by examining the effects
of changes in the dilution rate $D$ and the initial nutrient
concentration $N(0)$.  These are logical choices as the authors
 \citet{bec051} associated the chaotic states with the value of the
dilution rate.  Our previous work led us to choose the initial
nutrient concentration as well, as non-steady access to resources can
affect the overall behavior of the system.  We calculated the maximum
Lyapunov exponent for the system using $D$ in the range $[0.3, 2]$ and
$N(0)$ in the range $[0.1,1]$, and found that it was positive across
this entire range.

Figure \ref{becks-pc} contains parallel coordinates plots generated
using ranges of values for $D$, $N(0)$, and the corresponding maximum
Lyapunov exponent of each sample point.  By including maximum Lyapunov
exponent in these plots, we see that there is a small ``island'' of
sample points that behaves qualitatively differently from the rest;
when we isolate these points by increasing the lower bound on maximum
Lyapunov exponent, we see that they come from a very specific range of
values in terms of $D$ and $N(0)$.

\begin{figure}[h!]
\pdfgraphic{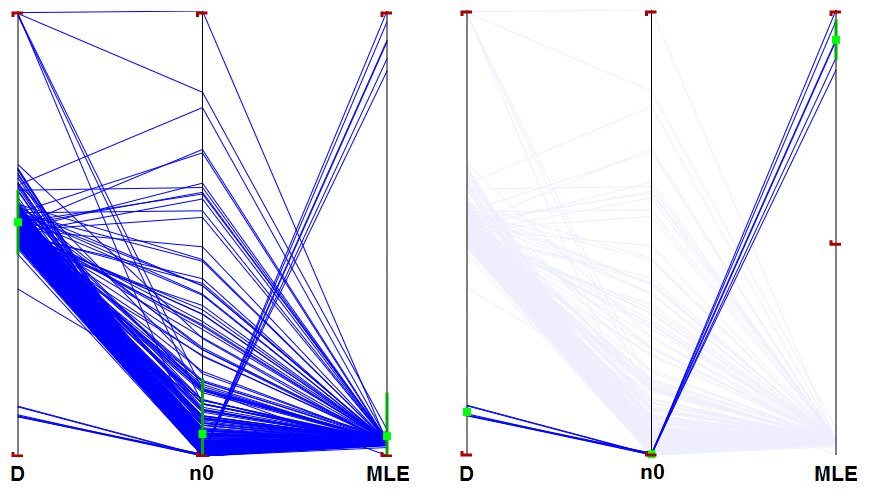}{0.5}{becks-pc}{Parallel coordinates plots of $D$, $N(0)$, and
maximum Lyapunov exponent (MLE).  The plot on the right shows the result of dragging
the lower bound on MLE upward to isolate only sample points with large MLE.}
\end{figure}

We also generated bifurcation graphs associated with this system of
equations using the dimensionalized equations and varying the dilution
rate $D$.  These are shown in Figure \ref{Becks-bifurcation}.  The
plots were generated using the MATLAB \verb=ode45= ODE integrator and
plotting tools.  The system was solved over the time frame
$t\in[0,7500]$ for each choice of $D$.  The plot contains the
solutions computed for the two prey species, $R$ and $C$, versus $D$
for $t\in[7000,7500]$ using 500 evenly spaced points.  We see evidence
of chaotic behavior for $D\in[0.2, 1.2]$.
\begin{figure}[h!]
\pdfgraphic{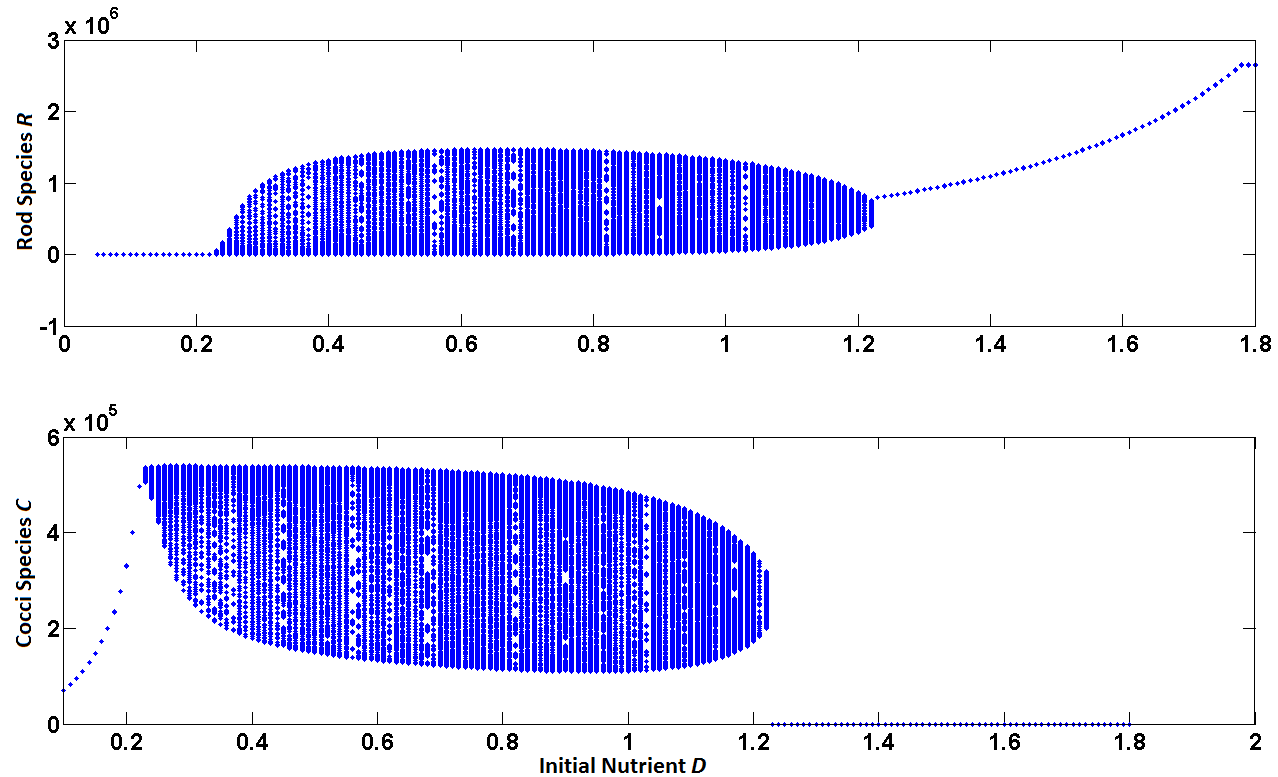}{0.4}{Becks-bifurcation}{Bifurcation plot of the ODE system showing the solutions for $R$ and $C$ versus $D$. Note the effect of the changing dilution rate.}
\end{figure}

\section{Conclusions}
\label{sec:conc}

We have provided a framework that allows researchers to efficiently
search over a design space including parameter values and initial
conditions and discover possible connections between values of these
constants and chaotic dynamics. The ability for such a search may
reveal chaotic behavior in systems not previously known to have
chaotic regime and reveal the existence of parameters and initial
conditions not previously known to yield chaotic behavior in studied
systems. We described a software framework for those wishing to
analyze systems for chaotic dynamics without the need for rigorous
mathematical analysis. One of its main strengths is 
allowing for studying system behavior that is virtually impossible to
observe in laboratory environment, thus making it useful to
experimentalists. The parallel coordinate plots further provide
insight into understanding system dynamics for a wide variety of
parameter changes.

Our algorithms can be further optimized. One obvious direction is to
investigate further methods of parallelizing the environment, which
would significantly speed up the computations.

\section*{Acknowledgements}

The authors express gratitude to Dr. Oleg Yordanov for his guidance and help with appropriate system rescaling; and to Dr. Fred Molz for conversations about chaos and Monoid kinetics, as well as initial development of the model equations. 

\bibliographystyle{plainnat}

\end{document}